# Approximate Joint Diagonalization and Geometric Mean of Symmetric Positive Definite Matrices

Marco Congedo[1]*, Bijan Afsari[2], Alexandre Barachant[1], Maher Moakher[3]

**1** GIPSA-lab, CNRS and Grenoble University, Grenoble, France, **2** Center for Imaging Science, Johns Hopkins University, Baltimore, Maryland, United States of America, **3** LAMSIN, National Engineering School, Tunis, Tunisia

* marco.congedo@gmail.com

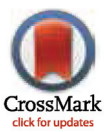

**Data Availability Statement:** All data required to replicate our study is available within the manuscript.

**Funding:** Author MC in an investigator of the European project ERC-2012-AdG-320684-CHESS and for this research has been partially supported by it. The support consisted in the reimbursement of expenses related to a scientific mission (visit to author MM) and the payment of publication fees. No other funder has supported this research besides the employers of the authors. The funders had no role in study design, data collection and analysis, decision to publish, or preparation of the manuscript.

## Abstract

We explore the connection between two problems that have arisen independently in the signal processing and related fields: the estimation of the geometric mean of a set of symmetric positive definite (SPD) matrices and their approximate joint diagonalization (AJD). Today there is a considerable interest in estimating the geometric mean of a SPD matrix set in the manifold of SPD matrices endowed with the Fisher information metric. The resulting mean has several important invariance properties and has proven very useful in diverse engineering applications such as biomedical and image data processing. While for two SPD matrices the mean has an algebraic closed form solution, for a set of more than two SPD matrices it can only be estimated by iterative algorithms. However, none of the existing iterative algorithms feature at the same time fast convergence, low computational complexity per iteration and guarantee of convergence. For this reason, recently other definitions of geometric mean based on symmetric divergence measures, such as the Bhattacharyya divergence, have been considered. The resulting means, although possibly useful in practice, do not satisfy all desirable invariance properties. In this paper we consider geometric means of covariance matrices estimated on high-dimensional time-series, assuming that the data is generated according to an instantaneous mixing model, which is very common in signal processing. We show that in these circumstances we can approximate the Fisher information geometric mean by employing an efficient AJD algorithm. Our approximation is in general much closer to the Fisher information geometric mean as compared to its competitors and verifies many invariance properties. Furthermore, convergence is guaranteed, the computational complexity is low and the convergence rate is quadratic. The accuracy of this new geometric mean approximation is demonstrated by means of simulations.

## Introduction

The study of distance measures between symmetric positive definite (SPD) matrices and the definition of the center of mass for a number of them has recently grown very fast, driven by practical problems in radar data processing, image processing, computer vision, shape analysis, medical imaging (especially diffusion MRI and Brain-Computer Interface), sensor networks, elasticity, mechanics, numerical analysis and machine learning (e.g., [1–10]). Interestingly, in this endeavor disparate perspectives from matrix analysis, operator theory, differential geometry, probability and numerical analysis have yielded converging results. In fact, we arrive at the same formalism and we end up with the same Riemannian metric from a pure *differential geometric* point of view [3, 5, 11–16], or from an *information geometric* point of view, assuming the multivariate Normal distribution of the data and adopting the *Fisher Information metric* [17, 18], dating back to the seminal works of Rao [19] and Amari [20].

The set of SPD matrices with a given dimension forms a smooth manifold, which we name the *SPD manifold*. As we will discuss below, there are various ways of defining natural geometries on the SPD manifold. If one is given a set of SPD matrices it may be useful to define a *center of mass* or *geometric mean* for the set taking into account the specific geometry of the manifold. One expects such a mean to be a more appropriate representative of the set as compared to other means that are not specific to the SPD manifold, such as the usual arithmetic and harmonic ones. The geometric mean is defined based on a chosen *metric* (distance) and several metrics have been proposed. These include the aforementioned Fisher information metric, the log-Euclidean metric [5, 21, 22] and the Bhattacharyya divergence [23, 24], also named *S*-divergence [25, 26], which turns out to be a specific instance of the $\alpha$-divergence [27]. The geometric mean based on the Fisher information metric, hereafter referred to as the FI (Fisher Information) mean, satisfies a number of desirable invariance properties, including congruence invariance, self-duality, joint homogeneity and the determinant identity. This is not always the case for geometric means based on other metrics, thus the FI mean is a fundamental object. Whereas for two SPD matrices the geometric mean has a straightforward definition, this is not the case for a set composed of more than two matrices [2]. For its estimation, regardless of the definition, one has to resort to either geometrical (constructive) procedures [28–30] or iterative optimization algorithms (see [30] for a comparison). In this work we do not consider constructive procedures because, in general, they do not satisfy all desirable properties of a geometric mean or, if they do, their computational cost becomes prohibitive as the number of matrices in the set increases [30]. Among iterative optimization algorithms for estimating the FI mean the most widely used is a simple gradient descent algorithm in the SPD manifold [10]. This algorithm features a moderate computational complexity per iteration and linear convergence rate. However, convergence itself is not guaranteed without choosing a small enough step-size. Specifically, if the algorithm convergences it does so to the only critical point, which is the global minimum, but it may not converge. Due to the convexity of the cost-function, in order to guarantee convergence the step-size must be adjusted as a function of the radius of a ball containing the data points [31] or an Armijo step-size search must be carried out at each iteration [30], thus the overall computational cost to ensure convergence is of concern. Recent attempts to apply conjugate gradient optimization and second order methods based on the exact or an approximate Hessian (such as trust-region and BFGS) improve the convergence rate, as expected. In [32] the authors show that for 3 x 3 matrices (e.g., diffusion tensors) a Newton algorithm based on explicit Hessian computations outperforms the gradient descent algorithm unless the radius of the ball is small. However, this advantage is completely nullified by the increased complexity per iteration as the dimension of the matrices in the set increases; extensive simulations performed by [30] have shown that with matrix dimension as little as ten the gradient descent approach is overall advantageous over all second order alternatives. Moreover, for the second order methods, convergence conditions are more restrictive. The recently proposed majorization-minimization algorithm of [33] guarantees convergence with linear convergence rate but high complexity per iteration, burdening its usefulness in practice for matrices of large dimension and/or large matrix sets. In light of this situation, the search for an efficient algorithm for estimating the geometric mean in large data set problems is currently a very active field.

In this article we introduce a new approximation to the FI mean springing from the study of the relation between the geometric mean of a set of SPD matrices and its approximate joint diagonalization [34–36]. We show that the invariance properties of this approximation derive from the invariance properties of the AJD algorithm employed. For instance, we obtain an approximation satisfying congruence invariance, joint homogeneity and the determinant identity, that is, all important properties of the FI geometric mean except self-duality, using the AJD algorithm developed in [36]. Using this AJD algorithm convergence is guaranteed, the computational complexity per iteration is low and the convergence rate is quadratic when the signal-to-noise ratio is favorable. As such, it offers an interesting alternative to existing iterative algorithms. Moreover, an on-line implementation is straightforward, allowing a fast on-line estimation of the geometric mean of an incoming stream of SPD matrices. We mention here that the approximate joint diagonalization (AJD) problem has originated in a completely different context and in the previous literature is completely unrelated to the problem of estimating the geometric mean. In fact, a solution to the AJD problem has arisen in statistics almost 30 years ago to solve the common principal component problem [37]. In the literature on signal processing it has appeared about 20 years ago [34] as an instrument to solve a very wide family of blind source separation problems, including the well-known independent component analysis [38]. Nonetheless, the connection between the geometric mean of a SPD matrix set and their AJD has remained unexplored so far.

In the following sections we introduce the notation and nomenclature. Then we review some concepts from Riemannian geometry and relevant metrics used to define a geometric mean. Then we establish the connection between the geometric mean of a SPD matrix set and their approximate joint diagonalizer. We introduce our approximation and we study its properties. In the result section we illustrate the accuracy of our approximation by means of simulations. Finally, we briefly discuss an on-line implementation and we conclude.

## Notation and Nomenclature

In the following we will indicate matrices by upper case italic characters ($A$), vectors, integer indices, random variables by lower case italic characters ($a$) and constants by upper case Roman characters (A). A set of objects will be enclosed in curly brackets such as $n \in \{1, \ldots, N\}$. Whenever possible, the same symbol will be used for indices and for the upper bound for summation and products, thus $\sum_n a_n$ will always stand short for $\sum_{n=1}^{N} a_n$. We will denote by $\text{tr}(\cdot)$, $|\cdot|$, $(\cdot)^T$, and $\|\cdot\|_F$ the trace of a matrix, its determinant, its transpose and its Frobenius norm, respectively. The operator $\text{diag}(\cdot)$ returns the diagonal part of its matrix argument. The identity matrix is denoted by $I$ and the matrix of zeros by $0$. $S$ will denote a symmetric matrix, $C$ a symmetric positive definite matrix (SDP) and $D$, $\Delta$ will be reserved for diagonal matrices. A set of K SPD matrices will be indicated by $\{C_1, \ldots, C_K\}$ or shortly as $\{C_k\}$. An asymmetric divergence from SPD matrix $C_2$ to SPD matrix $C_1$ will be denoted $\delta(C_1 \leftarrow C_2)$, whereas a symmetric distance or divergence between two SPD matrices will be denoted $\delta(C_1 \leftrightarrow C_2)$. The lambda symbol, as in $\lambda_n(A)$ will be reserved for the $n^{\text{th}}$ eigenvalue of matrix $A$. For the sake of brevity, notations of the kind $\ln^2\lambda_n(A)$ shall be read $(\ln\lambda_n(A))^2$. We will make extensive use of symmetric functions of eigenvalues for SDP matrices. For a symmetric matrix $S$ and SPD matrix $C$ these functions have general form $Uf(W)U^T$, where $U$ is the orthogonal matrix holding in its columns the eigenvectors of $S$ or $C$ and $W$ is the diagonal matrix holding the corresponding eigenvalues, to which the function applies element-wise. In particular, we will employ the following functions: inverse $C^{-1} = UW^{-1}U^T$, symmetric square root $C^{½} = UW^{½}U^T$, symmetric square root inverse $C^{-½} = UW^{-½}U^T$, logarithm $\ln(C) = U\ln(W)U^T$ and exponential $\exp(S) = U\exp(W)U^T$.

## Data Model

In many engineering applications we are confronted with multivariate observations obeying a linear instantaneous mixture generative model. For example, in electroencephalography (EEG) we observe N time-series of scalp electric potentials, typically sampled a few hundreds of times per second. Let us denote by $x(t) \in \Re^N$ the multivariate vector holding the observed data, in our example, the electric potentials recorded at the N scalp electrodes at time $t$. In EEG the sources of interest are equivalent electric dipoles formed by assemblies of locally synchronized pyramidal cells in the cortex. The current induced by these dipoles diffuses to the scalp sensors. Because of well-grounded physical and physiological reasons this process can be reasonably approximated as a linear instantaneous mixture [39], yielding generative model for the observed data

$$x(t) = As(t) + \eta(t), \qquad (1)$$

where $s(t) \in \Re^P$, $P \leq N$, holds the time series of the P cerebral sources to be estimated, $(t) \in \Re^N$ is a noise (model error) term, assumed uncorrelated to $s(t)$, and $A \in \Re^{N \times P}$ is the full column rank time-invariant *mixing matrix* determining the contribution of each source to each sensor, depending on the position and orientation of the electric dipoles. In the following we will assume for simplicity of exposition that P = N, i.e., we consider the estimation of as many sources as available sensors. Model (1) is quite general and appears in a very wide variety of engineering applications including speech, images, sensor array, geophysical and biomedical data processing. The well-known blind source separation (BSS) family of techniques, including independent component analysis [38], attempts to solve Equation (1) for $s(t)$, without assuming any knowledge on the mixing matrix $A$ (that is why it is named "blind"). The goal of BSS is to estimate the demixing matrix $B$ yielding source estimates $y(t)$ such as

$$y(t) = Bx(t) = BAs(t),$$

where we have ignored the noise term in (1) and where $A$ is the unknown true mixing matrix. An important family of BSS problems are solved by estimating a number of matrices holding second-order statistics (SOS) of multivariate observed measurements $x(t)$, e.g., in the case of EEG data, Fourier cospectral matrices estimated at different frequencies, covariance matrices estimated under different experimental conditions or covariance matrices estimated in different temporal windows, etc. [39]. According to model (1) such matrices, which are SPD, have theoretical form

$$C_k = AD_kA^T, \qquad (2)$$

where matrices $D_k \in \Re^{P \times P}$ here represent SOS matrices of the unknown P source processes and $k$ is the index for the K number of SOS matrices that are estimated. Under the assumption that the sources are uncorrelated (i.e., the matrices $D_k$ are diagonal), it can then be shown that the approximate joint diagonalization (AJD) of a set of K matrices generated under theoretical model (2) is indeed an estimation of the demixing matrix $B$ [40–42]. As per BSS theory, the sources can be estimated only up to an order and scaling indeterminacy, that is, we can at best find a matrix $B$ for which $BA$ ($A$ is the unknown true mixing process) approximates a matrix $P\Delta$, where $P$ is a permutation matrix and $\Delta$ an invertible diagonal matrix. This means that we can identify the *waveform* of the source processes, but their order and scaling (including sign switching) is arbitrary.

## Approximate Joint Diagonalization of SPD Matrices

The *joint diagonalizer* (JD) of two SPD matrices $C_1$ and $C_2$ is a matrix $B$ satisfying

$$\begin{cases} BC_1B^T = D_1 \\ BC_2B^T = D_2 \end{cases},$$
(3)

where $D_1$ and $D_2$ are diagonal matrices. The JD is not unique, since if $B$ is a JD of $C_1$ and $C_2$ so is $P\Delta B$, for any invertible diagonal matrix $\Delta$ and any permutation matrix $P$. The JD is obtained by the well-known generalized eigenvalue-eigenvector decomposition, as the matrix holding in the rows the eigenvectors of

$$C_1^{-1}C_2,$$
(4)

or the matrix holding in the rows the eigenvectors of

$$C_1^{-1/2}C_2C_1^{-1/2}$$
(5)

right-multiplied by $C_1^{-1/2}$, where indices 1 and 2 can be permuted in the above expressions. The JD matrix $B$ is orthogonal iff $C_1$ and $C_2$ commute in multiplication ([43], p. 160–165), which is not the case in general. Let $B = A^{-1}$ and $C_1, C_2$ be appropriate SOS matrices estimated from data $x(t)$, then conjugating by $A$ both sides of (3) we obtain the second order statistics of the generative model (2) (for the case K = 2) such as

$$\begin{cases} C_1 = AD_1A^T \\ C_2 = AD_2A^T. \end{cases}$$
(6)

In general, for a set of K>2 SPD matrices $\{C_1, \ldots, C_K\}$ there does not exist a matrix diagonalizing all of them. However, one may seek a matrix diagonalizing the set as much as possible, that is, we seek a matrix $B$ such that all products $BC_kB^T$, with $k \in \{1, \ldots, K\}$, are as diagonal as possible, according to some criterion. Such a problem is known as *approximate joint diagonalization* (AJD) and has been studied extensively in the signal processing community (e.g., [34], [36], [44]). As per the JD, the AJD matrix $B$ is not unique, since if $B$ is the AJD of set $\{C_1, \ldots, C_K\}$ so is $P\Delta B$, for any invertible diagonal matrix $\Delta$ and any permutation matrix $P$. As for the geometric mean, there is no closed-form expression for the AJD in the case K>2, so we proceed by specifying an AJD criterion and iteratively optimizing it. One such criterion specific for SPD matrices has been proposed by Pham [35, 36]; a consequence of the Hadamard inequality is that any SPD matrix $C$ verifies $|C| \leq |diag(C)|$, with equality iff $C$ is diagonal. Also, according to the Kullback-Leibler divergence $\delta_K(C \leftarrow D)$, the closest diagonal matrix $D$ to $C$ is $D = diag(C)$, for only in this case the divergence is zero [41]. The criterion proposed by Pham is then the sum of the Kullback-Leibler divergences of the input matrices to their diagonal form. Therefore we write a JD cost function as,

$$J(B|\{C_1, \cdots, C_K\}) = \sum_k \left[ \varsigma_k \left( \ln|diag(BC_kB^T)| - \ln|BC_kB^T| \right) \right],$$
(7)

where the $\varsigma_k$ are optional non-negative real numbers weighting the diagonalization effort with respect to the input matrices $C_k$. Besides being specific to SPD matrices, criterion ([7](#)) is interesting because it possesses an important invariance property, as stated in the following:

**Proposition 1.** Criterion ([7](#)) is invariant under positive rescaling of the input matrices $C_k$, thus the AJD matrix $B$ of $\{a_1 C_1, \ldots, a_K C_K\}$ is the same as the AJD matrix of $\{C_1, \ldots, C_K\}$ for any positive set $\{a_1, \ldots, a_K\}$.

The proof is trivial using well-known properties of the determinant and of the logarithm and will be omitted here. For AJD solutions according to a given criterion, we make use of the following:

**Definition 1.** The AJD of set $\{C_1, \ldots, C_K\}$ according to some AJD criterion is *well defined* (or is *essentially unique*) if any two global minimizers $B_1$ and $B_2$ of the criterion are in relation $B_1 = P\Delta B_2$, where $P$ is a permutation matrix and $\Delta$ is an invertible diagonal matrix.

For details on the essential uniqueness of AJD see [[40](#), [45](#), [46](#)]. This is a generic property in the sense that in the presence of noise the AJD is essentially unique (with very high probability). Now, consider again data model ([1](#)) and the theoretical generative model of SOS statistics of the observed data as per ([2](#)), where $A$ is the (unknown) mixing matrix and the $D_k$, with $k \in \{1, \ldots, K\}$, are the (unknown) diagonal matrices holding K SOS statistics of the source processes. Working with AJD we are interested in the situation where the noise term in ([1](#)) is small enough so that the AJD solution is well defined as per definition 1. We will then make use of a general property of such well-defined AJD solutions:

**Proposition 2.** For any invertible matrix $F$, if $B$ is an essentially unique AJD of set $\{FC_1 F^T, \ldots, FC_K F^T\}$, then $\Delta PBF$ is the AJD of set $\{C_1, \ldots, C_K\}$, where invertible diagonal matrix $\Delta$ and permutation matrix $P$ are the usual AJD scaling and permutation indeterminacy, respectively.

**Proof.** Saying that $B$ is an essentially unique AJD of set $\{FC_1 F^T, \ldots, FC_K F^T\}$ according to some AJD criterion implies that the set $\{BFC_1 F^T B^T, \ldots, BFC_K F^T B^T\}$ is a global minimizer of the AJD criterion employed. Thus, matrix $BF$ is, out of possible different scaling and permutation as per definition 1, a global minimizer of the same AJD criterion for the set $\{C_1, \ldots, C_K\}$.

Finally, we will need the following:

**Definition 2.** Let $B$, with inverse $A$, be an essentially unique AJD of set $\{C_1, \ldots, C_K\}$; an AJD criterion is said to verify the self-duality invariance if $P\Delta A^T$ is a well defined AJD of set $\{C_1^{-1}, \ldots, C_K^{-1}\}$ satisfying the same criterion, with invertible diagonal matrix $\Delta$ and permutation matrix $P$ the usual AJD scaling and permutation indeterminacy, respectively.

## The Riemannian Manifold of SPD Matrices

Geometrically, the Euclidean space of SPD matrices of dimension N x N can be considered as a ½N(N+1)-dimensional hyper cone ([Fig 1](#)). The usual vector (Euclidian) space of general square matrices is endowed with the metric $\langle S_1, S_2 \rangle = tr(S_1^T S_2)$ and associated (Frobenius) norm $\|S\|_F$. We will replace the convex pointed cone in the vector space of [Fig 1](#) with a regular manifold of non-positive curvature without boundaries, developing instead infinitely in all of its ½N(N+1) dimensions. In differential geometry, a N-dimensional smooth manifold is a topological space that is locally similar to the Euclidean space and has a globally defined differential structure. A smooth *Riemannian manifold* (or *Riemannian space*) M is a real smooth manifold equipped with an inner product on the tangent space $T_\Omega M$ defined at each point $\Omega$ that varies smoothly from point to point. The tangent space $T_\Omega M$ at point $\Omega$ is the Euclidean vector space containing the tangent vectors to all curves on M passing through $\Omega$ ([Fig 2](#)). In the SPD manifold endowed with the Fisher information metric for any two vectors $V_1$ and $V_2$ in

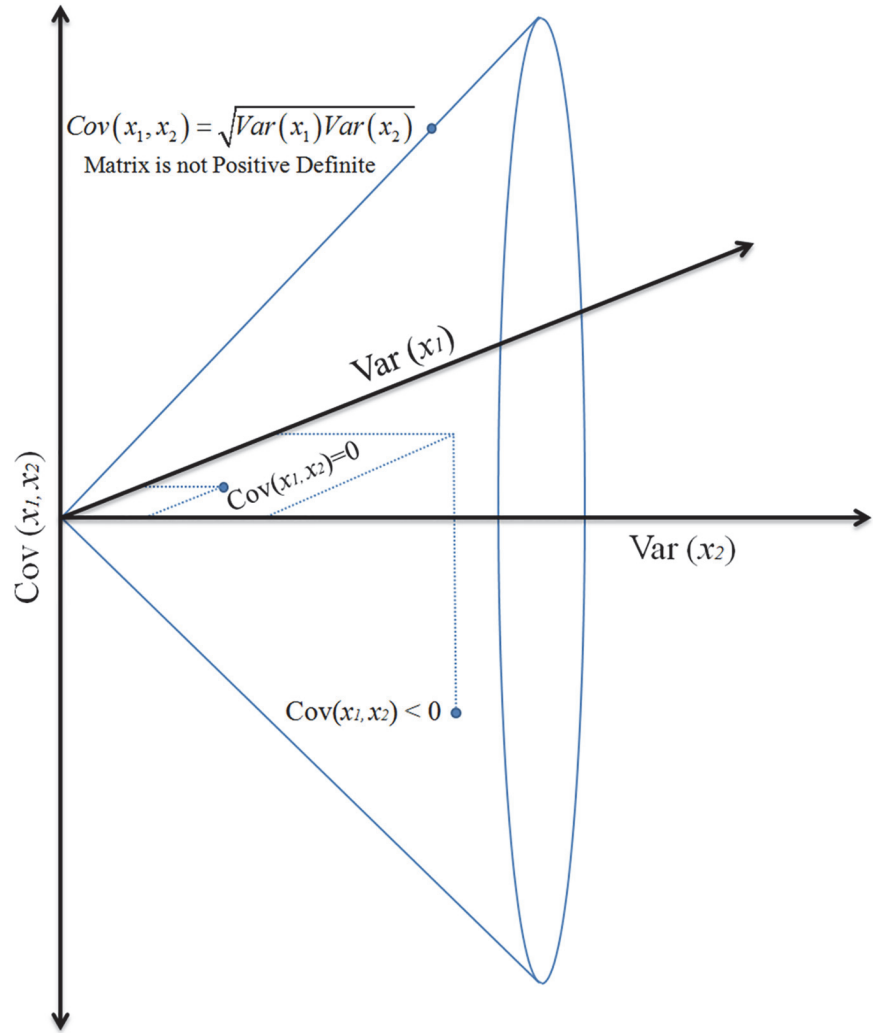

**Fig 1. Symmetric positive definite matrices, e.g. covariance matrices, are constrained by their symmetry, the strict positivity of the diagonal elements (variance) and the Cauchy-Schwarz inequalities bounding the absolute value of the off-diagonal elements:** $|Cov(x_i|x_j)| \leq [Var(x_j)Var(x_i)]^{1/2}$, **for all i,j∈{1,...,N}.** This topology is easily visualized in case of 2x2 matrices; any 2x2 covariance matrix can be seen as a point in 3D Euclidean space, with two coordinates given by the two variances (diagonal elements) and the third coordinate given by the covariance (either one of the off-diagonal element). By construction a covariance matrix must stay within the cone boundaries. As soon as the point touches the boundary of the cone, the matrix is no more positive definite.

the tangent space the inner product through point $\Omega$ is given by [2]

$$\langle V_1, V_2 \rangle_\Omega = tr(\Omega^{-1} V_1 \Omega^{-1} V_2).$$

**The Geodesic.** The Fisher information metric allows us to measure the length of curves in M and find the shortest curve between two points $\Omega$ and $\Phi$ on M. This is named the geodesic and is given by [2]

$$\gamma_\beta(\Omega \to \Phi) = \Omega^{1/2}(\Omega^{-1/2} \Phi \Omega^{-1/2})^\beta \Omega^{1/2}, \ \beta \in [0, 1], \tag{8}$$

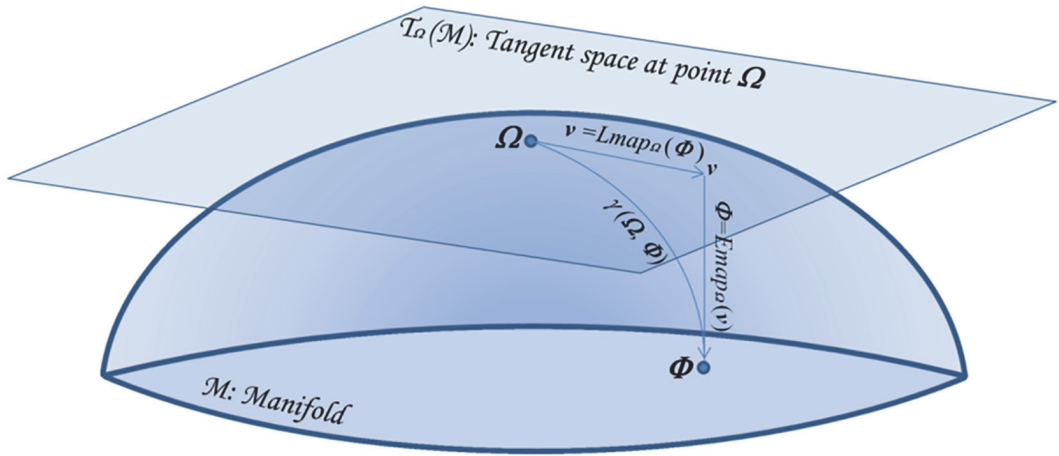

**Fig 2. The Manifold and the tangent space at a point.** Consider a point $\Omega$ on M and construct the tangent space $T_\Omega M$ on it. Now take a tangent vector $V$ departing from $\Omega$, which is our reference point. There exists one and only one geodesic on the manifold starting at $\Omega$ that corresponds to $V$; think at rolling the plane (tangent space) on the surface (manifold) in such a way that the vector always touches the surface. The end point on M is $\Phi$. We see that the geodesics on M through $\Omega$ are transformed into straight lines and the distances along all geodesics are preserved (this is true in the neighborhood of $\Omega$). (Rearranged from [17]).

doi:10.1371/journal.pone.0121423.g002

where $\beta$ is the arc length parameter. When $\beta = 0$ we are at $\Omega$, when $\beta = 1$ we are at $\Phi$ and when $\beta = 1/2$ we are at the *geometric mean* of the two points (Fig 2).

**The Exponential and Logarithmic Maps.** The exponential and logarithmic maps are shown graphically in Fig 2. The function that maps a vector $V \in T_\Omega M$ to the point $\Phi \in M$ following the geodesic starting at $\Omega$, is named the *exponential map* and denoted by $\Phi = \text{Emap}_\Omega(V)$. It is defined as

$$\Phi = \text{Emap}_\Omega(V) = \Omega^{1/2} \exp\left(\Omega^{-1/2} V \Omega^{-1/2}\right) \Omega^{1/2} \qquad (9)$$

The inverse operation is the function mapping the geodesic relying $\Omega$ to $\Phi$ back to the tangent vector $V \in T_\Omega M$. It is named the *logarithmic map* and denoted $V = \text{Lmap}_\Omega(\Phi)$. It is defined as

$$V = \text{Lmap}_\Omega(\Phi) = \Omega^{1/2} \ln\left(\Omega^{-1/2} \Phi \Omega^{-1/2}\right) \Omega^{1/2}. \qquad (10)$$

**The Metric (Distance).** Given two points $C_1$ and $C_2$ on the manifold $M$, their *Riemannian distance* based on the Fisher information metric is the length of the geodesic (8) connecting them. It is given by [3, 10, 16]

$$\delta_R(C_1 \leftrightarrow C_2) = \|\ln(C_2^{-1} C_1)\|_F = \sqrt{tr \ln^2(\Lambda)} = \sqrt{\sum_n \ln^2 \lambda_n}, \qquad (11)$$

where $\Lambda$ is the diagonal matrix holding the eigenvalues $\lambda_1, \ldots \lambda_N$ of either matrix (4) or (5). This distance has a remarkable number of properties, some of which are reported in Table 1 [14, 15, 25]. For more inequalities see [2, 25]. Associated to the chosen metric is also the *Riemannian norm*, defined as the Riemannian distance from an SPD matrix $C$ to the identity, that is, the Euclidean distance of its logarithm to the zero point:

$$\|C\|_R = \delta_R(I \leftrightarrow C) = \|\ln(C)\|_F = \sum_n \ln^2 \lambda_n(C). \qquad (12)$$

The Riemannian norm is zero only for the identity matrix, while the Frobenius norm is zero

**Table 1. Some important properties and inequalities of the Riemannian Fisher Information Distance.**

| Fundamental Properties of the Riemannian Fisher information distance |
| --- |
| (13) Positivity $\delta_R(\Omega \leftrightarrow \Phi) \geq 0$, with equality iff $\Omega = \Phi$ |
| (14) Symmetry $\delta_R(\Omega \leftrightarrow \Phi) = \delta_R(\Phi \leftrightarrow \Omega)$ |
| (15) Congruence-Invariance $\delta_R(\Omega \leftrightarrow \Phi) = \delta_R(B\Omega B^T \leftrightarrow B\Phi B^T)$, for any invertible $B$ |
| (16) Similarity-Invariance $\delta_R(\Omega \leftrightarrow \Phi) = \delta_R(B^{-1}\Omega B \leftrightarrow B^{-1}\Phi B)$, for any invertible $B$ |
| (17) Invariance under Inversion $\delta_R(\Omega \leftrightarrow \Phi) = \delta_R(\Omega^{-1} \leftrightarrow \Phi^{-1})$ |
| (18) Proportionality $\delta_R(\Omega \leftrightarrow \gamma_\beta(\Omega \to \Phi)) = \beta\delta_R(\Omega \leftrightarrow \Phi)$ |
| **Some inequalities of the Riemannian Fisher information distance** |
| (19) $\delta_R(\gamma_\beta(\Omega \to \Phi) \leftrightarrow \gamma_\beta(\Omega \to \Xi)) \leq \beta\delta_R(\Phi \leftrightarrow \Xi)$ |
| (20) $\delta_R(\Omega \leftrightarrow \Phi) \geq \|\ln\Omega - \ln\Phi\|_F$, with equality iff $\Omega$ and $\Phi$ commute |

doi:10.1371/journal.pone.0121423.t001

only for the null matrix. Either an eigenvalue smaller or greater than 1 increases the norm and the norm goes to infinity as any eigenvalues go to either infinity or zero. Importantly, because of the square of the log, an eigenvalue $\lambda$ increases the norm as much as an eigenvalue $1/\lambda$ does (see Fig 3), from which the invariance under inversion.

**The Geometric Mean of SPD Matrices: general considerations.** Given a set of SPD matrices $\{C_1, \dots, C_K\}$, in analogy with the arithmetic mean of random variables, a straightforward definition of the matrix arithmetic mean is

$$\mathcal{A}\{C_1, \cdots, C_K\} = \frac{1}{K}\sum_k C_k, \qquad (13)$$

and a straightforward definition of the harmonic mean is

$$\mathcal{H}\{C_1, \cdots, C_K\} = \left(\frac{1}{K}\sum_k C_k^{-1}\right)^{-1}.$$

On the other hand a straightforward definition of the geometric mean is far from obvious because the matrices in the set in general do not commute in multiplication. Researchers have postulated a number of desirable properties a mean should possess. Ten such properties are

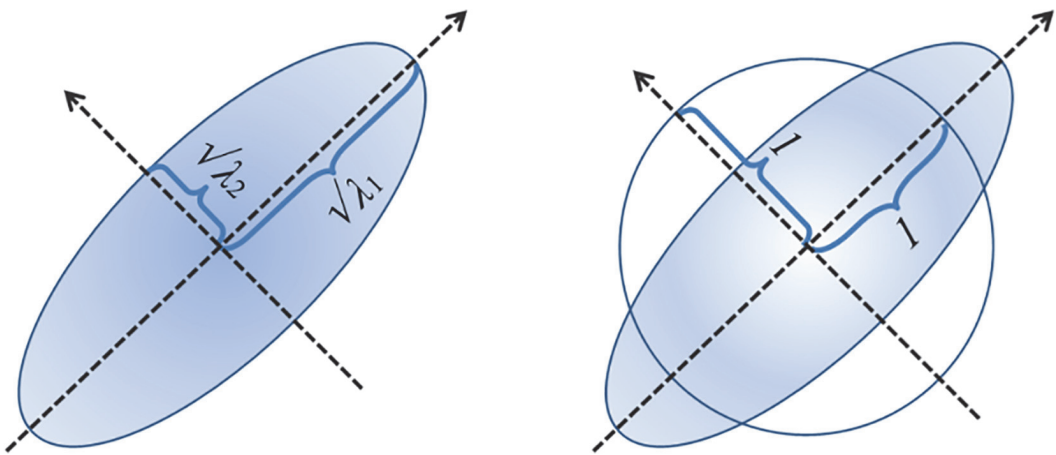

**Fig 3. The ellipsoids in the figure are isolines of constant density of bivariate Gaussian distributions.** The semiaxes are proportional to the square root of the eigenvalues of the covariance matrix. If we ask how far the ellipsoid is from the circle, which is the definition of the norm (12), we see that an eigenvalue = 2 contribute to the distance from the identity as much as an eigenvalue = 1/2 does, as one would expect, since the eigenvalues are squared quantities. Neither the sum nor the sum of the logarithm of the eigenvalue verify this property.

doi:10.1371/journal.pone.0121423.g003

**Table 2. Some important properties of the Fisher information metric geometric mean.**

| Properties of the geometric Mean |
| --- |
| (25) Invariance by Reordering: $\mathcal{G}_R\{C_1, \cdots, C_K\}$ is the same for any order of matrices in the set |
| (26) Congruence Invariance: $F(\mathcal{G}_R\{C_k\})F^T = \mathcal{G}_R\{FC_kF^T\}$, for any invertible $F$ |
| (27) Self-Duality: $\mathcal{G}_R^{-1}\{C_k\} = \mathcal{G}_R\{C_k^{-1}\}$ |
| (28) Joint Homogeneity: $\mathcal{G}_R\{a_1C_1, \cdots, a_KC_K\} = \left(\prod_k a_k\right)^{1/K}\mathcal{G}_R\{C_1, \cdots, C_K\}$, $a_k \geq 0$ |
| (29) Determinant Identity: $|\mathcal{G}_R\{C_1, \cdots, C_K\}| = \left(\prod_k |C_k|\right)^{1/K}$ |
| (30) if all matrices $C_k$ pair-wise commute then |
| $\mathcal{G}_R\{C_1, \cdots, C_K\} = \left(\prod_k C_k\right)^{1/K} = \exp\left(1/K \sum_k ln(C_k)\right)$ |
| (31) $\delta_R^2(\mathcal{G}_R\{C_k\} \leftrightarrow \Omega) \leq \sum_k 1/K [\delta_R^2(C_k \leftrightarrow \mathcal{G}_R\{C_k\}) - \delta_R^2(C_k \leftrightarrow \Omega)]$, for any PSD $\Omega$ |

doi:10.1371/journal.pone.0121423.t002

known in the literature as the ALM properties, from the seminal paper in [11]. The concept of Fréchet means and the ensuing variational approach are very useful in this context: in a univariate context the arithmetic mean minimizes the sum of the squared Euclidean distances to K given values while the geometric mean minimizes the sum of the squared hyperbolic distances to K given (positive) values. In analogy, we define the (least-squares) Riemannian geometric mean $G\{C_1, \ldots, C_k\}$ of K SPD matrices $C_k$ such as the matrix satisfying [13, 15].

$$\underset{\mathcal{G}\{C_1, \cdots, C_K\}}{argmin} \sum_k \delta^2(\mathcal{G}\{C_1, \cdots, C_K\} \leftrightarrow C_k), \qquad (14)$$

where $\delta^2(\cdot \leftrightarrow \cdot)$ is an appropriate squared distance. In words, the Riemannian geometric mean is the matrix minimizing the sum of the squared distances of all elements of the set to itself. Using the Fisher information (FI) distance (11) such mean, which we name the *FI mean* and denote as $\mathcal{G}_R\{C_1, \cdots, C_K\}$ or, shortly, $\mathcal{G}_R\{C_k\}$, features all ALM properties. The FI mean is the *unique* SPD *geometric mean* satisfying non-linear matrix equation [15]

$$\sum_k ln(\mathrm{G}_R^{-1}\{C_1, \cdots, C_K\}C_k) = 0 \qquad (15)$$

or, equivalently,

$$\sum_k \ln\left(\mathcal{G}_R^{-1/2}\{C_1, \cdots, C_K\}C_k\mathcal{G}_R^{-1/2}\{C_1, \cdots, C_K\}\right) = 0, \qquad (16)$$

where, in line with the notation used in this article, $\mathcal{G}^s\{\cdot\}$ stands short for $(\mathcal{G}\{\cdot\})^s$, for any power s. We have listed a number of properties of the FI mean along with some related inequalities in Table 2 (see also [2, 25, 47, 48]). Note that in the literature this least-squares geometric mean is sometimes referred to as the barycenter, the Riemannian center of mass, the Fréchet mean, the Cartan mean or the Karcher mean (although these definitions in general are not equivalent and attention must be paid to the specific context, see e.g., [31]).

**The Geometric Mean and the Joint Diagonalization of Two SPD Matrices.** Given two points $C_1$ and $C_2$ on the manifold M, the Geometric Mean of them, indicated in the literature by $C_1\#C_2$, has several equivalent closed-form expressions, such as [2, 12, 25, 28, 47]

$$C_1\#C_2 = C_1^{1/2}\left(C_1^{-1/2}C_2C_1^{-1/2}\right)^{1/2}C_1^{1/2} = C_1(C_1^{-1}C_2)^{1/2} = (C_2C_1^{-1})^{1/2}C_1 \qquad (17)$$

And

$$C_1\#C_2 = C_1^{1/2}\exp\left(1/2 \ln\left(C_1^{-1/2}C_2C_1^{-1/2}\right)\right)C_1^{1/2} = C_1\exp\left(1/2 \ln(C_1^{-1}C_2)\right). \qquad (18)$$

In the above the indices 1 and 2 can be switched to obtain as many more expressions. The geometric mean of two SPD matrices is indeed the midpoint of the geodesic in (8) and turns out to be the unique solution of a quadratic Ricatti equation [2, 48], yielding

$$(C_1 \# C_2)C_2^{-1}(C_1 \# C_2) = C_1 \text{ and } (C_1 \ C_2)C_1^{-1}(C_1 \# C_2) = C_2.$$

Our investigation has started with the following:

**Proposition 3.** The FI geometric mean of two SPD matrices can be expressed uniquely in terms of their joint diagonalizer; let $B$ be the JD (3) of matrices $\{C_1, C_2\}$ such that $BC_1B^T = D_1$ and $BC_2B^T = D_2$ and let $A = B^{-1}$, then the geometric mean is

$$(C_1 \# C_2) = A(D_1 D_2)^{1/2} A^T, \tag{19}$$

for any JD solution, regardless the permutation and scaling ambiguity.

**Proof.** Since diagonal matrices commute in multiplication, using (3) and properties (50) and (54) of the geometric mean we can write $B(C_1 \# C_2)B^T = D_1 \# D_2 = (D_1 D_2)^{1/2}$. Conjugating both sides by $A = B^{-1}$ we obtain (19).

**Remark.** A consequence of the scaling indeterminacy of the JD and (19) is that we can always chose $B$ such that

$$BC_1B^T BC_2B^T = I, \tag{20}$$

in which case $A$ is a square root ([49], p. 205–211) of the geometric mean, i.e.,

$$(C_1 \# C_2) = AA^T, \text{ given (20) true.} \tag{21}$$

**The Geometric Mean of a SPD Matrix Set.** Given a set $\{C_1, \ldots, C_K\} = \{C_K\}$ of $K > 2$ SPD matrices the point $G_R\{C_k\}$ satisfying (23)–(24) has no closed form solution. The research of algorithms for estimating the FI geometric mean or a reasonable approximation for the case $K > 2$ is currently a very active field [3, 24, 28–30, 48, 50]. In this work we consider two iterative algorithms for estimating the FI mean. The first is the aforementioned gradient descent algorithm [10]. In order to estimate $G_R\{C_k\}$, we initialize $M$ with the arithmetic mean or any other smart guess. Then we apply the iterations given by

$$M \leftarrow M^{1/2} \exp\left[\upsilon/K \sum_k \ln\left(M^{-1/2} C_k M^{-1/2}\right)\right] M^{-1/2}, \tag{22}$$

until convergence. Here $\upsilon$ is the step size, typically fixed at 1.0. Notice that this algorithm iteratively maps the points in the tangent space through the current estimation of the mean (10), computes the arithmetic mean in the tangent space (where the arithmetic mean makes sense) and maps back the updated mean estimation on the manifold (9), until convergence (Fig 4). Note also that the global minimum attained by gradient descent (22) satisfies

$$\sum_k \ln\left(M^{-1/2} C_k M^{-1/2}\right) = 0, \tag{23}$$

upon which the iterate does not move anymore. We have already reported that these iterations have linear convergence rate, but convergence itself is not guaranteed. In order to minimize the occurrence of divergence and also to speed-up convergence, while avoiding computationally expensive searches of the optimal step size, we use in the present work the following version with heuristical decrease of the step-size over iterations:

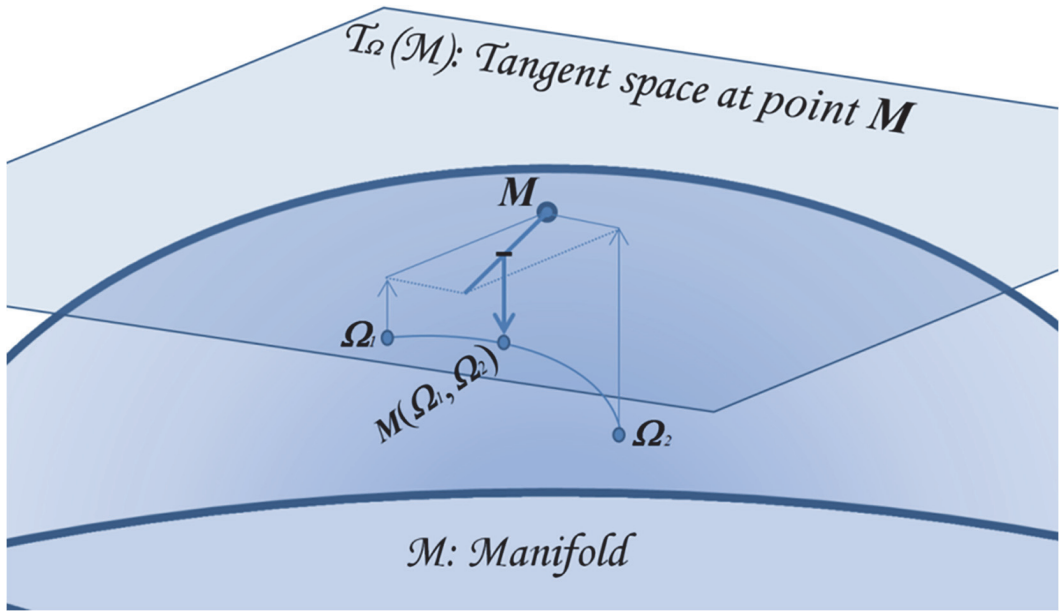

**Fig 4. Zoom on the manifold as it is represented in Fig 2.** Consider two points $\Omega_1$ and $\Omega_2$ on M and construct the tangent space $T_G$M through the current estimation of their mean $M$, initialized as the arithmetic mean. At each iteration, the algorithm maps the points on the tangent space, computes the mean vector and maps back the point on the manifold. At each iteration the estimation of the mean is updated, thus the point of transition into the tangent space changes, until convergence, that is, until this transition point will not change anymore, coinciding with the geometric mean, that is, satisfying (23).

## GM-GD Algorithm

Initialize M

Set ε equal to a suitable machine floating point precision number (e.g., $10^{-9}$ for double precision),

$\upsilon = 1$, τ equal to the highest real number of the machine.

Repeat

$$h \leftarrow \left\| \upsilon/K \sum_k \ln\left(M^{-1/2} C_k M^{-1/2}\right) \right\|_F$$

*If $h < \tau$*

*then* $M \leftarrow M^{1/2} \exp\left[\upsilon/K \sum_k \ln\left(M^{-1/2} C_k M^{-1/2}\right)\right] M^{1/2}$, $\upsilon \leftarrow 0.95\upsilon$, $\tau \leftarrow h$

*else* $\upsilon \leftarrow 0.5\upsilon$

*Until* $\left(\left\| 1/K \sum_k \ln\left(M^{-1/2} C_k M^{-1/2}\right) \right\|_F < \epsilon\right)$ OR ($\upsilon < \epsilon$)

The second algorithm we consider is the aforementioned majorization-minimization algorithm recently proposed in [33]. Since convergence is guaranteed, it is used here as a benchmark for accuracy:

## GM-MM Algorithm

Initialize M

Repeat

$$\left|\begin{aligned}
&\text{For } k \text{ from 1 to K do } \Psi_k \leftarrow \ln\left(C_k^{-1/2} M C_k^{-1/2}\right)\\
&\Phi_1 \leftarrow \sum_k \left(C_k^{1/2}\left((\Psi_k^2 + I)^{1/2} - \Psi_k\right)C_k^{-1/2}M^{-1}\right)\\
&\Phi_2 \leftarrow \sum_k \left(C_k^{-1/2}\left((\Psi_k^2 + I)^{1/2} + \Psi_k\right)C_k^{1/2}M^{-1}\right)\\
&M \leftarrow \Phi_1^{1/2}\left(\Phi_1^{1/2}\Phi_2\Phi_1^{1/2}\right)^{-1/2}\Phi_1^{1/2}
\end{aligned}\right. \tag{24}$$

Until Convergence

**Alternative Metrics and Related Geometric Means.**    Recently it has been proposed to use the least-squares (Fréchet) geometric mean ([14]) based on the log-Euclidean distance [21, 22, 51] and the Bhattacharyya divergence [23, 24], also named $S$-divergence [25, 26], which turns out to be a specific instance of the $\alpha$-divergence setting $\alpha = 0$ [27]. These studies suggest that these alternative definitions of geometric mean give results similar to the FI mean in practical problems, while their computational complexity is lower.

The *Log-Euclidean distance* is

$$\delta_L(C_1 \leftrightarrow C_2) = \|\ln(C_1) - \ln(C_2)\|_F, \tag{25}$$

which is the straightforward generalization to matrices of the scalar hyperbolic distance and, from property (48), equals the FI distance if $C_1$ and $C_2$ commute in multiplication. The interest of this distance is that the geometric (Fréchet) mean ([14]) in this case has closed-form solution given by

$$\mathcal{G}_L\{C_1, \cdots, C_K\} = \exp\left(1/K \sum_k \ln(C_k)\right). \tag{26}$$

This is a direct generalization of the geometric mean of positive scalars and, again, from property (48), it equals the FI mean if all matrices in the set pair-wise commute. The computation of (26) requires fewer flops than a single iteration of the gradient descent algorithm ([22]), so this constitutes by far the most efficient geometric mean estimation we know. The log-Euclidean mean (26) possesses the following properties: invariance by reordering, self-duality, joint homogeneity and determinant identity [21–22]. It does not possess the congruence invariance, however, it is invariant by rotations (congruence by an orthogonal matrix). Also, important for the ensuing derivations, whereas the determinant of the log-Euclidean mean is the same as the determinant of the FI mean (since both verify the determinant identity), the trace of the log-Euclidean mean is larger than the trace of the FI mean, unless all matrices in the set pair-wise commute, in which case as we have seen, the two means coincide, hence they have the same trace and determinant [21]. Because of this trace-increasing property, the log-Euclidean mean is in general farther away from the FI mean as compared to competitors such as constructive geometric means [28]. Moreover, in general the log-Euclidean mean differs from the FI mean even for the case K = 2.

In [24] the author has shown that the *Bhattacharyya divergence* (also named log-det divergence and S-divergence) between two SPD matrices $C_1, C_2 \in \Re^{N \cdot N}$ behaves similarly to the FI distance if the matrices are close to each other. It is symmetric as a distance, however, it does not verify the triangle inequality. This divergence is

$$\delta_B^2(C_1 \leftrightarrow C_2) = \ln\frac{|(C_1+C_2)/2|}{1/2\ln|C_1 C_2|} = \ln\frac{|\mathcal{A}\{C_1,C_2\}|}{|\mathcal{G}\{C_1,C_2\}|} = \ln\prod_n 1/2\left(\sqrt{\lambda_n}+\frac{1}{\sqrt{\lambda_n}}\right), \quad (27)$$

where the arithmetic mean $\mathcal{A}\{C_1,C_2\}$ is defined in(21), the geometric mean $\mathcal{G}\{C_1,C_2\} = C_1 \ C_2$ in(32)–(33) and the $\lambda_n$ are, again, the N eigenvalues of either matrix (4) or (5). Successively, in [26] the author has shown that the *square-root* of the Bhattacharyya divergence is a distance metric satisfying the triangle inequality. Both the Bhattacharyya divergence and distance are invariant under inversion and under congruence transformation [24, 26]. The geometric mean (14) of a set of matrices based on divergence (27) is the solution to nonlinear matrix equation [26]

$$2/\text{K}\sum_k(\mathcal{G}_B\{C_1,\cdots,C_\text{K}\}+C_k)^{-1} = \mathcal{G}_B^{-1}\{C_1,\cdots,C_\text{K}\}. \quad (28)$$

In order to estimate $\mathcal{G}_B\{C_1,\cdots,C_\text{K}\}$ we initialize $M$ and apply iterations [24]

$$M \leftarrow \text{K}\left[\sum_k\left(\frac{C_k+M}{2}\right)^{-1}\right]^{-1}, \quad (29)$$

which global minimum has been shown to satisfy (28) [26]. This mean possesses the following invariance properties: invariance by reordering, congruence invariance, self-duality and joint homogeneity. However, it does not satisfy the determinant identity. As a consequence, in general *both* the determinant and the trace of the Bhattacharyya mean differ from those of the FI mean. However, the Bhattacharyya mean of K = 2 matrices coincides with the FI mean [27].

**The Geometric Mean of a SPD Matrix Set by AJD.**   Let $H$ be an invertible matrix of which the inverse is $G$. Using the congruence invariance of the geometric mean (50) and conjugating both sides by $G$ we obtain

$$\mathcal{G}_R\{C_1,\cdots,C_\text{K}\} = G\mathcal{G}_R\{HC_1H^T,\cdots,HC_\text{K}H^T\}G^T.$$

Our idea is to reach an approximation of $\mathcal{G}_R\{C_1,\cdots,C_\text{K}\}$ by approximating $\mathcal{G}_R\{HC_1H^T,\cdots,HC_\text{K}H^T\}$ in the expression above. Particularly, if the matrices $HC_kH^T$ are nearly diagonal, then they nearly commute in multiplication, hence we can employ property (54) to approximate the geometric mean $\mathcal{G}_R\{C_1,\cdots,C_\text{K}\}$ by expression

$$\mathcal{G}_A'' = A\mathcal{G}_L\{BC_1B^T,\cdots,BC_\text{K}B^T\}A^T = A\exp\left(1/\text{K}\sum_k\ln(BC_kB^T)\right)A^T, \quad (30)$$

where $\mathcal{G}_L\{BC_1B^T,\cdots,BC_\text{K}B^T\}$ is the log-Euclidean mean introduced in (41), $B$ is a well-defined AJD matrix of the set $\{C_1,\ldots,C_\text{K}\}$ (definition 1) chosen so as to minimize criterion (7), and $A$ is its inverse. Because of the scaling ($\Delta$) and permutation ($P$) ambiguities of AJD solutions, (30) actually defines an infinite family of admissible means. Before we start to study the specific instance of the family we are interested in (the closest to the FI mean in general), let us observe that if all products $BC_kB^T$ are exactly diagonal then the family of means defined by (30) collapses on a single point in the manifold, which is indeed the FI geometric mean, for any AJD matrix with form $P\Delta B$. This is the case when the data are generated according to model (1) and the noise term therein is null, in which case the coincidence of the two means is a consequence of property (54) or, regardless the data model, if K = 2, in which case (30) reduces to (19). The same is true also whenever

$$1/\text{K}\sum_k\ln(BC_kB^T) = 0, \quad (31)$$

in which case the matrix exponential of (30) is the identity matrix and the family (30) collapses to the unique point $AA^T$ just as for the case K = 2 in (21). Interestingly, the family of means (30) includes the log-Euclidean mean, in that setting $A = B^{-1} = I$ we obtain (26), but also the FI mean, in the sense that setting $A = M^{1/2}$ and $B = A^{-1} = M^{-1/2}$ we obtain the global minimum of the FI mean (23). This is the point on the manifold we sought to approximate by approximating condition (31). First of all, we show that the AJD permutation ambiguity is not of concern here. We have the following

**Proposition 4.** The family of means given by (30) is invariant with respect to the AJD permutation indeterminacy $P$, for any invertible AJD solution $B$ with inverse $A$ and any invertible diagonal scaling matrix $\Delta$.

**Proof.** Taking into account the $P$ and $\Delta$ indeterminacies, since for any permutation matrix $P^T = P^{-1}$ and for any diagonal matrix $\Delta = \Delta^T$, the family of solutions (30) has full form

$$A\Delta^{-1}P^{-1}\exp\left(1/K \sum_k \ln(P\Delta BC_k B^T \Delta P^{-1})\right)P\Delta^{-1}A^T.$$

If f (S) is a function of the eigenvalues of S (see section "Notation and Nomenclature") and $P$ is an arbitrary permutation matrix, then $P^{-1}f(PSP^{-1})P = f(S)$, thus, the family of solutions (30) actually takes form

$$\mathcal{G}'_A\{C_1, \cdots, C_K\} = A\Delta^{-1}\exp\left(1/K \sum_k \ln(\Delta BC_k B^T \Delta)\right)\Delta^{-1}A^T \tag{32}$$

for any invertible $B = A^{-1}$ and any invertible diagonal $\Delta$.

Then, note that the family of means defined by (32) is SPD, which is a consequence of the fact that the exponential function of a symmetric matrix is always SPD and that both $B$ and $\Delta$ are full rank. It also verifies the invariance under reordering (49), which is inherited from the invariance under reordering of AJD algorithms and of the log-Euclidean mean. Furthermore, we have the following

**Proposition 5.** The family of means defined by (32) verifies the determinant identity(53), for any invertible matrix $B = A^{-1}$ and any invertible diagonal $\Delta$.

**Proof.** We need to prove that for any invertible matrix $B = A^{-1}$ and any invertible scaling matrix $\Delta$

$$|\mathcal{G}'_A\{C_1, \cdots, C_K\}| = |\mathcal{G}_R\{C_1, \cdots, C_K\}| = (\prod_k |C_k|)^{1/K}.$$

Using (32) we write the left-hand side of the above equation such as $|A||\Delta^{-1}||\exp\left(1/K \sum_k \ln(\Delta BC_k B^T \Delta)\right)||\Delta^{-1}||A^T|$. Since the log-Euclidean mean possesses the determinant identity, this is

$$|A||\Delta^{-1}|(\prod_k |\Delta BC_k B^T \Delta|)^{1/K}|\Delta^{-1}||A^T|.$$

Developing the products of determinants and since $|A| = |A^{-1}|^{-1} = |B|^{-1}$, we obtain the desired result such as

$$|A||\Delta^{-1}||\Delta||B|(\prod_k |C_k|)^{1/K}|B^T||\Delta||\Delta^{-1}||A^T| = (\prod_k |C_k|)^{1/K}.$$

**Proposition 6.** If $B$ is the AJD solution of (7), with inverse $A$, the family of means defined by (32) verifies the joint homogeneity property (52), for any invertible diagonal $\Delta$.

**Proof.** We need to prove that

$$\mathcal{G}'_A\{a_1 C_1, \cdots, a_K C_K\} = \left(\prod_k a_k\right)^{1/K} \mathcal{G}'_A\{C_1, \cdots, C_K\}$$

The result follows immediately from the invariance under rescaling of criterion ([7]) (Proposition 1), as

$$A\Delta^{-1}\left[\exp\left(1/K \sum_k \ln(\Delta B a_k C_k B^T \Delta)\right)\right]\Delta^{-1}A^T$$
$$= \left(\prod_k a_k\right)^{1/K} A\Delta^{-1}\left[\exp\left(1/K \sum_k \ln(\Delta B C_k B^T \Delta)\right)\right]\Delta^{-1}A^T.$$

So far we have considered the properties of the whole family of means with general form, ([32]) for which we have solved the permutation AJD indeterminacy (Proposition 4). We now seek the member of the family better approximating the FI geometric mean given a matrix $P\Delta B$ that performs as the AJD of the set $\{C_1, \ldots, C_K\}$. This involves choosing a specific scaling $\Delta$. As we have seen, if in ([32])

$$\exp\left(1/K \sum_k \ln(\Delta B C_k B^T \Delta)\right) = I, \qquad (33)$$

then ([21]) is true also for K>2, in which case the FI geometric mean is given by $AA^T$, where $A$ is the inverse of $\Delta B$ in ([32]) and ([30]) is a stationary point of ([22]). Such condition is well approximated if the left-hand side of ([33]) is nearly diagonal and all its diagonal elements are equal, that is, in our case, by scaling the rows of $B$ such that

$$\text{diag}\left[\exp\left(1/K \sum_k \ln(B C_k B^T)\right)\right] = \alpha I, \quad \text{for an arbitrary } \alpha > 0. \qquad (34)$$

The uniqueness of this solution, regardless the choice of $\alpha$, is demonstrated by the following

**Proposition 7.** Let $B = A^{-1}$ be an invertible AJD of the set $\{C_1, \ldots, C_K\}$. Scaling $B$ so as to satisfy ([34]), the mean ([32]) is unique and invariant with respect to $\alpha$.

**Proof.** Once matrix $B$ satisfies ([34]), let us consider a further fixed scaling such us $\Delta = \omega I$, with $\omega > 0$. Substituting this fixed scaling matrix in the family of means ([32]) we have

$$A\omega^{-1}I\exp\left(1/K \sum_k \ln(\omega IB C_k B^T I\omega)\right)I\omega^{-1}A^T,$$

but this is

$$A\omega^{-1}(\omega^{2K})^{1/K}\omega^{-1}\exp\left(1/K \sum_k \ln(B C_k B^T)\right)A^T = A\exp\left(1/K \sum_k \ln(B C_k B^T)\right)A^T,$$

showing that the resulting mean does not change. Thus, given an AJD matrix $B$, we obtain our sought unique member of the family as

$$\mathcal{G}_A\{C_1, \cdots, C_K\} = A\exp\left(1/K \sum_k \ln(B C_k B^T)\right)A^T, \text{ given ([34]) true for any } \alpha > 0, \qquad (35)$$

where $A$ is the inverse of the matrix $B$ scaled so as to satisfy ([34]).

Without loss of generality, hereafter we will choose $\alpha = 1$. With this choice the AJD matrix $B$ is an approximate "whitening" matrix for the sought mean, since we have

$$\text{diag}[B\mathcal{G}_A(C_1, \cdots, C_K)B^T] = \text{diag}\left[BA\left[\exp\left(1/K \sum_k \ln(B C_k B^T)\right)A^T B^T\right]\right] = I.$$

In order to satisfy ([34]), notice that any change in $\Delta$ changes the log-Euclidean mean therein since the log-Euclidean mean is not invariant under congruent transformations, thus there is no closed-form solution to match condition (34) given an arbitray AJD matrix solution. However, we easily obtain the desired scaling of $B$ by means of an iterative procedure (see below). We name the resulting approximation to the geometric mean (35) the ALE mean, where ALE is the acronym of *AJD-based log-Euclidean Mean*. The algorithm is as follows:

## ALE Mean Algorithm

Let $B$ be an invertible AJD of set $\{C_1, \ldots, C_K\}$ minimizing criterion (7).

Set $\epsilon$ equal to a suitable machine floating point precision number

Repeat

$$\left| \begin{aligned} &\Delta \leftarrow \text{diag}\left[\exp\left(1/K \sum_k \ln(BC_kB^T)\right)\right] \\ &B \leftarrow \Delta^{-1/2}B \end{aligned} \right. \tag{36}$$

Until $^1/_N \delta_R(I \leftrightarrow \Delta) \leq \epsilon$;

$$A \leftarrow B^{-1}$$

$$\mathcal{G}_A\{C_1, \cdots, C_K\} = A\exp\left(1/K\sum_k \ln(BC_kB^T)\right)A^T$$

Note that instead of the average FI distance $^1/_N \delta_R(I \leftrightarrow \Delta)$ any suitable distance of $\Delta$ from the identity matrix can be used as well in order to terminate the Repeat-Until loop here above. Note also that if instead of (36) we perform iterations

$$B \leftarrow \left[\exp\left(\upsilon/K \sum_k \ln(BC_kB^T)\right)\right]^{-1/2}B, \tag{37}$$

(which is the same as (36) without the "diag" operator) the algorithm converges to an inverse square root of the FI geometric mean, i.e., upon convergence (and, in this case, if the algorithm converges)

$$\mathcal{G}_G\{C_1, \cdots, C_K\} = (B^TB)^{-1}. \tag{38}$$

In fact, (37) is the gradient descent equivalent to (22) converging to the inverse square root of the geometric mean instead of on the geometric mean itself. Equivalently, one may iterate

$$A \leftarrow A\left[\exp\left(\upsilon/K \sum_k \ln(BC_kB^T)\right)\right]^{1/2}, \tag{39}$$

converging to the square root of the FI geometric mean (as before, $A$ here is a square root of the FI mean). In the case of iterations (37) or (39), however, there is no point to perform a previous AJD ($B$ in (37) and $A$ in (39) can be initialized, for instance, by the inverse square root and the square root of the arithmetic mean, respectively) and we encounter the same convergence dependency on the step-size $\upsilon$ as for the gradient descent (22). Our approximation based on (36) instead surely converges without any step-size to be searched, as convergence of the AJD algorithm is ensured without using a step-size [36] and iterations (36) imply only a scaling of $B$ and also surely converge. Besides providing a unique solution, our ALE mean satisfies a very important property: whereas the log-Euclidean mean is not invariant under congruent transformation, the ALE mean is. This is demonstrated by the following:

**Proposition 8.** The ALE mean (35) satisfies the invariance under congruent transformation.

**Proof.** We need to show that for any invertible matrix $F$

$$F\mathcal{G}_A\{C_1, \cdots, C_K\}F^T = \mathcal{G}_A\{FC_1F^T, \cdots, FC_KF^T\}.$$

Let $B_1$ a well defined AJD of set $\{C_1, \ldots, C_K\}$ with inverse $A_1$ and $B_2$ a well-defined AJD of set $\{FC_1F^T, \ldots, FC_KF^T\}$ with inverse $A_2$, both satisfying condition (35) for their respective set. The expression above is then

$$FA_1\exp\left(1/\text{K}\sum_k\ln(B_1C_kB_1^T)\right)A_1^TF^T = A_2\exp\left(1/\text{K}\sum_k\ln(B_2FC_kF^TB_2^T)\right)A_2^T$$

Because of Proposition 2, if matrix $B_1$ approximately diagonalizes the set $\{C_1, \ldots, C_K\}$ so does matrix $B_2F$ according to the same criterion and if they both satisfy (35) for the set $\{C_1, \ldots, C_K\}$ they are equal out of a permutation indeterminacy that we can ignore because of proposition 4. As a consequence $A_1 = (B_2F)^{-1} = F^{-1}A_2$ and thus $A_2 = FA_1$. Making the substitutions we obtain

$$FA_1\exp\left(1/\text{K}\sum_k\ln(B_1C_kB_1^T)\right)A_1^TF^T = FA_1\exp\left(1/\text{K}\sum_k\ln(B_1C_kB_1^T)\right)A_1^TF^T.$$

**Proposition 9.** The ALE mean verifies the self-duality property (51) if the AJD solution $B$, with inverse $A$, verifies the self-duality property of Definition 2.

**Proof.** Self-duality of the ALE mean is verified if

$$\mathcal{G}_A\{C_1, \cdots, C_K\} = \mathcal{G}_A^{-1}\{C_1^{-1}, \cdots, C_K^{-1}\}.$$

Using definition 2 we have

$$\mathcal{G}_A^{-1}\{C_1^{-1}, \cdots, C_K^{-1}\} = \left[B^T\exp\left(1/\text{K}\sum_k\ln(A^TC_k^{-1}A)\right)B\right]^{-1},$$

and computing the inverse of the right-hand side

$$\mathcal{G}_A^{-1}\{C_1^{-1}, \cdots, C_K^{-1}\} = A\left[\exp\left(1/\text{K}\sum_k\ln(A^TC_k^{-1}A)\right)\right]^{-1}A^T.$$

Using $\ln(S^{-1}) = -\ln(S)$ and then $(\exp(-S))^{-1} = \exp(S)$ we have

$$\mathcal{G}_A^{-1}\{C_1^{-1}, \cdots, C_K^{-1}\} = A\exp\left(\frac{1}{K}\sum_k\ln(BC_kB^T)\right)A^T = \mathcal{G}_A\{C_1, \cdots, C_K\}.$$

Note that the AJD matrix satisfying criterion (7) verifies the self-duality of definition 2 only if K = 2 or if all products $B^TC_kB$ are exactly diagonal, otherwise it verifies it only approximately. Therefore, in general our ALE mean estimation (35) verifies self-duality (51) only approximately. Interestingly, the ALE mean would verify the self-duality property as well if the AJD cost function makes use of a diagonality criterion based on the Riemannian distance (11) instead of the Kullback-Leibler divergence as in (7). However, the search of such an AJD algorithm has proven elusive so far. In conclusion, we have shown that the ALE mean verifies several important properties satisfied by the FI Riemannian mean (invariance under reordering and under congruence transformation, joint homogeneity and determinant equality). The self-duality property is satisfied approximately in general and exactly in some special cases. Next we will study the performance of the ALE mean by means of simulations.

## Results

As mentioned in the introduction, the estimation of the FI geometric mean of a set of SPD matrices has proven useful in a number of diverse engineering applications [1–10]. For instance, working with brain-computer interfaces based on electroencephalography (EEG) we have found it very useful for defining prototypical brain states and classifying single trials according to the distance of the trials to different prototypes. Since EEG data is generated under the instantaneous linear mixing model ([1])[39], in this section we explore the possible advantage of estimating the geometric mean using the ALE estimation. In order to perform simulations we generate sets of SPD matrices according to model ([1]) (see ([2]), but also see the generation of simulated matrices in [28, 33]): a set of K matrices is generated as

$$C_k = 10(A_{True}D_kA_{True}^T + \Sigma_k),\qquad(40)$$

where $A_{True} \in \Re^{N \cdot N}$ (true mixing matrix) has entries randomly drawn from a standard Gaussian distribution, $D_k$ is a diagonal matrix with entries randomly drawn from a squared standard Gaussian distribution and bounded below by $10^{-4}$ and $\Sigma_k$ is a noise matrix obtained as $(1/N)QQ^T$ where $Q \in \Re^{N \cdot N}$ is a matrix with entries randomly drawn from a Gaussian distribution with mean zero and standard deviation $\sigma$ (noise level).

[Fig 5](#) shows typical convergence patterns for the GD (Gradient Descent with decreasing step-size), MM (Majorization-Minimization), Bha (Bhattacharyya) and Pham's AJD algorithm (Pham's AJD algorithm is the first part of the ALE algorithm) for K = 100 matrices of dimension N = 10, with three different noise levels ($\sigma$ = 0.01, 0.1, 1). These three noise levels correspond approximately to typical low, medium and high noise situations according to AJD standards, that is to say, the diagonalization that can be obtained on a set generated with $\sigma$ = 1 would be considered very bad for practical AJD purposes and the AJD matrix $B$ in this case cannot be assumed well-defined. We can appreciate the typical superlinear convergence rate of Pham's AJD algorithm in low and medium noise situations. For high noise situations the convergence rate becomes slower. The opposite happens for the GD, MM and Bha algorithm, which convergence rate becomes faster as the noise level increases. For considering the overall computational efficacy of these algorithms we should also consider the complexity per iteration; it is the lowest for Pham's AJD algorithm, followed in the order by the Bha, GD and MM algorithms. In fact, the log-Euclidean mean ([26]) necessitates the computation of K eigenvalue-eigenvector decompositions only once. On the other hand, Bha ([29]) involves K+1 Cholesky decompositions at each iteration, GD ([22]) involves K eigenvalue-eigenvector decompositions at each iteration and MM ([24]) involves 2K+2 eigenvalue-eigenvector decompositions at each iteration. Iterations ([36]) in the second part of the ALE algorithm converge reliably in a few iterations requiring K eigenvalue-eigenvector decompositions per iteration. Overall then, the ALE mean based on Pham's AJD is advantageous when Pham's algorithm convergence rate is quadratic. This happens when the noise in data model ([1]) is small enough.

[Fig 6](#) shows the typical relation observed between the trace and the determinant of the FI geometric mean estimated by algorithm GD and MM (they converge to the same point), LE, Bha and ALE. As we can see the log-Euclidean mean has always the same determinant as the FI mean, but larger trace. The Bhattacharyya mean has always both different trace and different determinant. The ALE mean has determinant and trace only very slightly different, regardless the noise level. As the noise increases the bias of the LE and Bha mean tends to decrease.

[Fig 7](#) shows the FI distance ([11]) between the LE, Bha, ALE geometric mean and the FI mean estimated by the GD and MM algorithms (they converge to the same point, which is used here as a benchmark) for K = 100 matrices of dimension N = 10, with three different noise level ($\sigma$ = 0.01, 0.1, 1). The distance is plotted against the condition number of the true mixing matrix in our generative model ([40]). The ALE mean is consistently closer to the FI mean in all cases and, surprisingly, in absolute terms it is pretty close even in the high-noise case ($\sigma$ = 1). Also, as suggested by [Fig 6](#), the estimation of LE and Bha approaches the FI mean as the noise increase. On

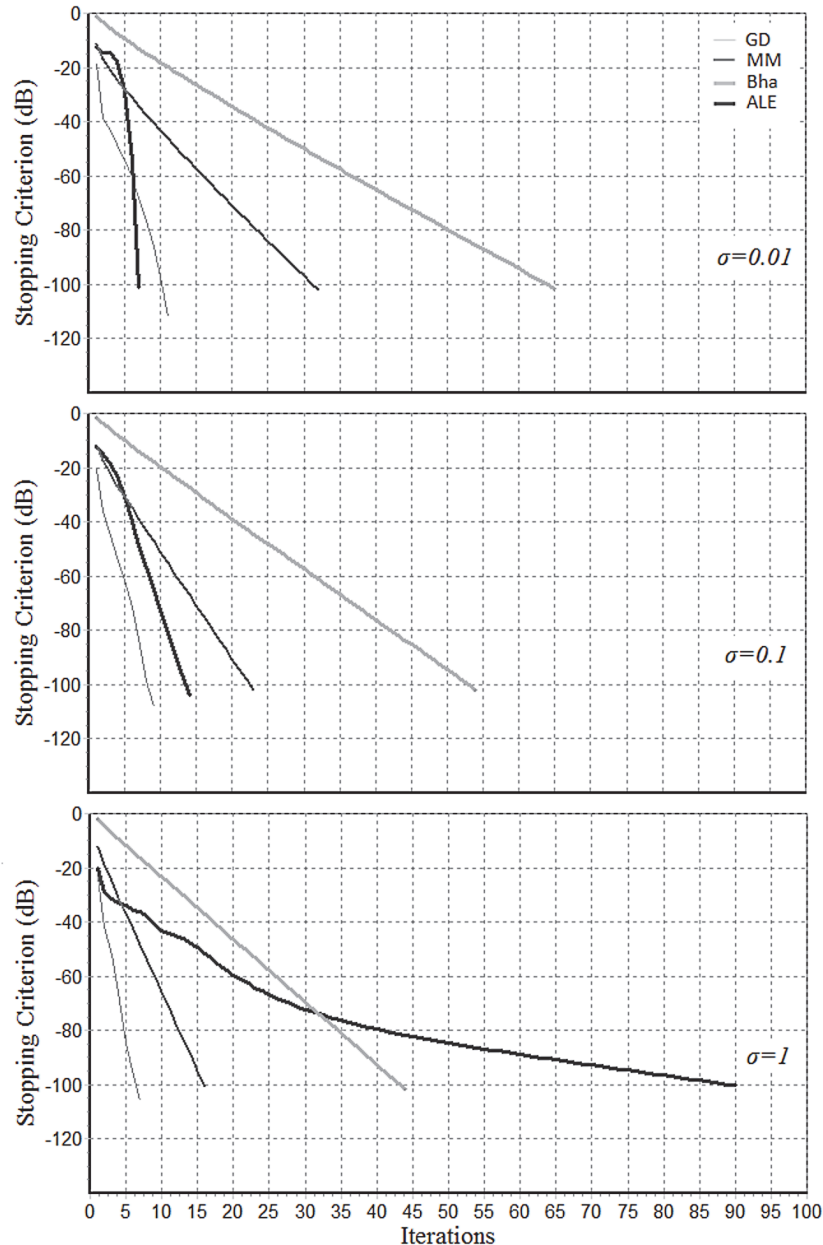

**Fig 5. Typical convergence behavior of algorithms GM, MM, Bha and ALE (first part of ALE algorithm, corresponding to AJD estimation by Pham's algorithm).** The simulations are done generating data according to (40) for three noise levels ($\sigma$ = 0.01, 0.1, 1). Each algorithm was stopped when its own stopping criterion became smaller than -100dB.

doi:10.1371/journal.pone.0121423.g005

the other hand, the conditioning number appears to play a role only for the low and moderate noise cases ($\sigma$ = 0.01, 0.1).

## Conclusions and Discussion

In this paper we explored the relationship between the approximate joint diagonalization of a SPD matrix set and its geometric mean. After appropriate scaling, the inverse of the joint diagonalizer of two SPD matrices is a square root of their geometric mean. For the general case of a

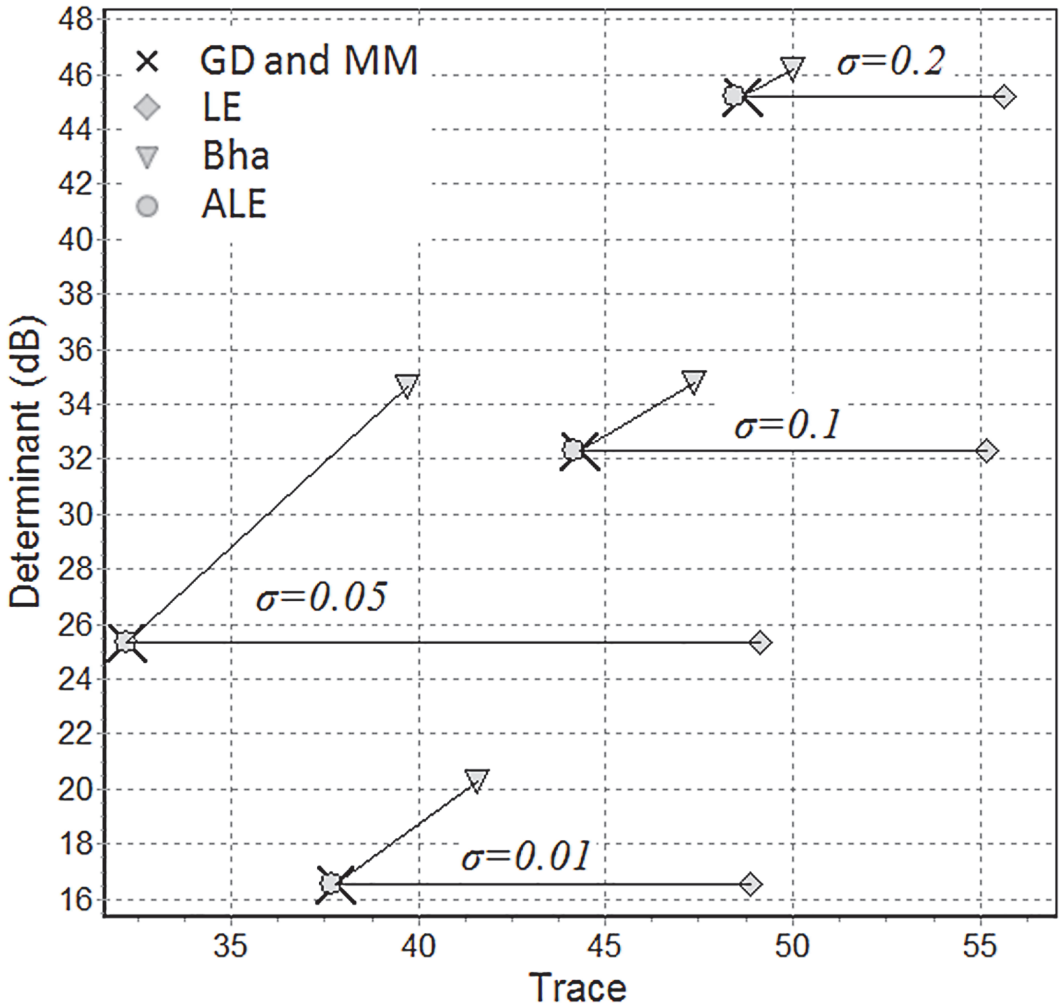

**Fig 6. Typical relation between the trace (x-axis) and the determinant (y-axis, in dB) of the FI geometric mean as estimated with algorithm GD and MM (cross), and the means estimated by LE (diamond), Bha (triangle) and ALE (disk) for four noise levels (σ = 0.01, 0.05, 0.1 0.2), N = 10 and K = 100.** The arrows link the FI estimation with the corresponding LE and Bha estimations.

SPD matrix set comprised of more than two matrices, we have studied a family of geometric means that includes the geometric mean according to the Fisher information metric (the FI mean) and the log-Euclidean mean. Then we have introduced a specific instance of this family, which is computationally attractive and does not require a search for the optimal step size. We have showed that it approximates the FI geometric mean much better than the log-Euclidean mean. Indeed, this mean, named the ALE mean, can be conceived as an improved version of the log-Euclidean mean, in that i) it satisfies the congruence invariance and ii) similar to the log-Euclidean mean it has the same determinant as the FI mean, but has much smaller trace, thus its trace is much closer to the trace of the FI mean. The ALE mean can be computed by running an AJD algorithm followed by a specific scaling obtained by a simple iterative procedure. The AJD algorithm developed by Pham [36] is particularly adapted for this purpose, as its convergence is guaranteed and features nice invariance properties, which translate into a number of invariance properties for the ALE mean. For this algorithm the convergence rate is quadratic when the set can be nearly diagonalized, that is, when the data is generated according to model (1) and the noise term therein is small enough. In such a situation our ALE mean is clearly advantageous over its competitors both in terms of stability (guarantee of convergence) and computational complexity. Also, for matrices of big dimension Pham's AJD algorithm can be easily parallelized, since it proceeds in a Jacobi-like fashion updating pair-wise the vectors of the matrix [36]. An on-line implementation is also straightforward and very efficient [41]: one creates a circular buffer of matrices and once an incoming matrix enters the buffer one or two iterations of the algorithm suffice to update the AJD solution and one or two iterations (36) suffice to update the ALE mean altogether. By applying appropriate weighting to the matrices in the buffer and/or regulating the buffer size one decides how fast the ALE mean should adapt to the incoming data. Our simulations have confirmed our analysis and have shown that the ALE mean can be employed even if the noise term in data generation model is high, i.e., even if the matrices are very far apart from each other in terms of their Riemannian distance; however in this case the computational advantage of the ALE mean vanishes. We conclude that the ALE mean is the convenient choice when the set can be diagonalized pretty well; otherwise the gradient descent algorithm is computationally advantageous, even if searching for the optimal step-size at each iteration.

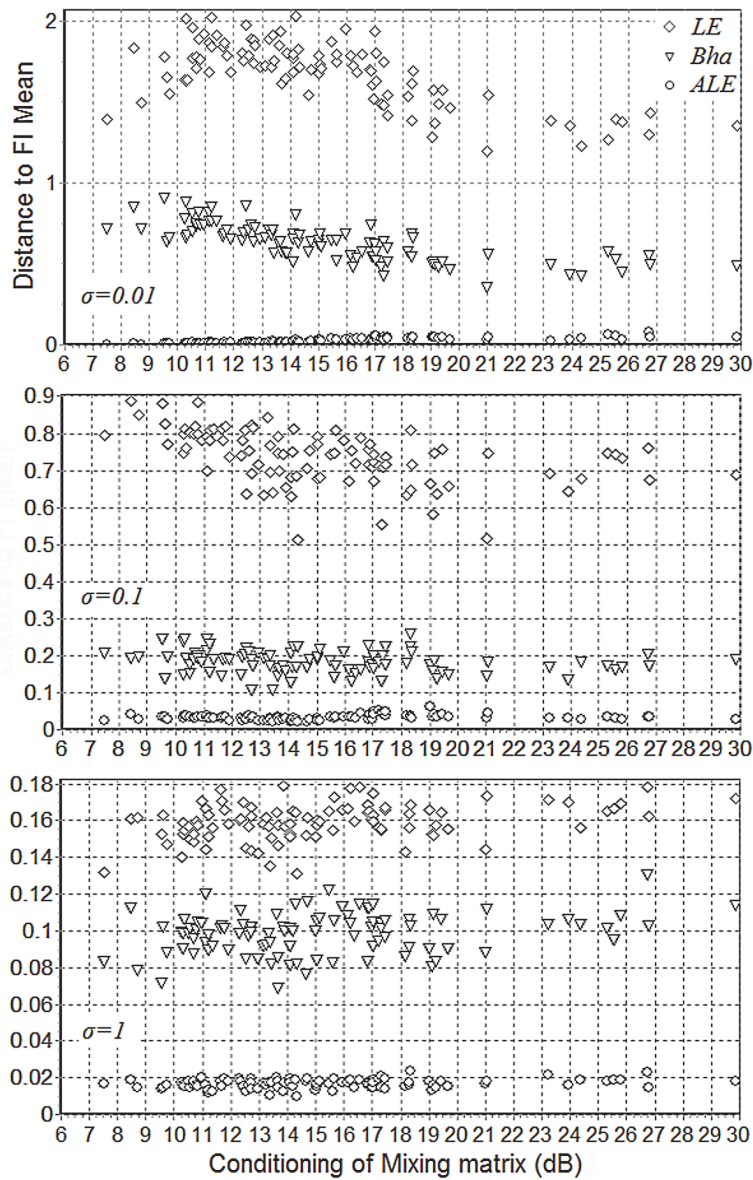

**Fig 7. Distance of the LE (diamonds), Bha (triangles) and ALE (disks) geometric mean to the FI geometric mean** estimated by the GD and MM algorithm (they converge to the same point), for three noise levels (σ = 0.01, 0.1, 1) and variable condition numbers of the true mixing matrix in (40). The simulations were repeated 100 times, with N = 10 and K = 100. Notice the different scales on the y-axes.

## Author Contributions

Conceived and designed the experiments: MC BA AB MM. Performed the experiments: MC AB. Analyzed the data: MC AB. Contributed reagents/materials/analysis tools: MC AB. Wrote the paper: MC BA AB MM.